\documentclass[a4paper,10pt]{article}
\usepackage{amssymb}
\usepackage{latexsym}

\begin{document}

\begin{center}
\LARGE { Moduli Spaces of Instantons on Noncommutative 4-Manifolds II} \\
\vspace{2cm}
\Large { Noriaki HAYAKAWA and Hiroshi TAKAI} \\
\vspace{5mm}
\Large { Department of Mathematics, \\
Tokyo Metropolitan University} \\
\end{center}

\vspace{3cm}

\begin{center}
\Large  Abstract 
\end{center}

\large  \quad Studied are moduli spaces of self dual connections on noncommutative 4-manifolds, especially deformation quantization of compact toric Riemannian 4-manifolds. Then such moduli spaces of irreducible modules associated with highest weights of compact connected semisimple Lie groups are smooth manifolds with dimension determined by their weights. This may be viewed as a generalization 
of Atiyah-Hitchin-Singer's classical result, Landi-Suijlekom's noncommutative 4-sphere case and the recent result of the second authors.

\newpage
\Large {\bf{\S1.~Introduction}}\large \quad Among many important topics of 
Yang-Mills theory, Atiyah-Hitchin-Singer [AHS] showed that the moduli spaces of  irreducible self-dual connections of principal buldles over compact spin 
4-manifolds are smooth manifolds with their dimension determined in terms of 
topological invariants related to ambient bundles. Recent development of 
noncommutative geometry serves many powerful devices toward determining a 
noncommutative analog of their result cited above. For example, Connes-Landi[CL] showed an existence of an isospectral deformation of compact spin Riemannian 
manifolds whose isometry groups have at least 2-torus subgroup. Related with 
this results, Connes-Devois$\cdot$Violette [CV] observed that a deformation 
quantization of manifolds could be viewed as the fixed point algebra of the set of all smooth functions from manifolds to noncommutative tori under certain 
action of tori. Suceedingly, Landi-Suijlekom [LS] computed concerning the 
noncommutative 4-sphere, the index of the Dirac operator on certain spin modules as well as the dimension of the instanton moduli space of an irreducible module. Recently, the second author [T] verified that given noncommutative principal 
bundles based on compact spin toric Riemannian 4-manifolds and compact connected semisimple Lie groups and their associated noncommutative smooth sections with highest weights, their instanton moduli spaces are locally smooth manifolds with their dimension determined by their weights and the ambient bundles and as a 
corollary if the spaces of noncommutative anti-self-dual harmonic 2-forms vanish for all their self-dual connections, they are smooth manifolds with more 
explicit dimensions.\\
   In this paper, we show that given a principal bundle based on compact toric 
Riemannian 4-manifold and a compact connected Lie group and its associated 
noncommutative smooth sections with a highest weight, its instanton moduli space is a smooth manifolds with its dimension determined by the highest weight and 
the ambient bundles. This result is used to compute the concrete example 
associated with the Hopf bundle over the complex projective 2-space. \\

\Large{\bf{\S2.~Noncommutative Yang-Mills Theory}} \\
\large In this section, we review briefly both commutative Yang-Mills theory due to Atiyah-Hitchin-Singer [AHS] and the noncommutative Yang-Mills theory due to Landi-Suijlekom [LS]. Let $M$ be a compact oriented toric Riemannian 4-manifold, $G$ a compact connected Lie group, and let $P$ a principal $G$-bundle over $M$. Suppose there exist a smooth action $\varphi$ from 2-torus $T^2$ into the isometry group $Iso(M,g)$ of $M$ with respect to its Riemannian metric $g$. By Connes-Devois$\cdot$Violette [CV], let $M_{\theta}$ be a deformation quantization of $M$ along $\theta$ as Frechet *-algebra. Actually, there exists a quantization 
map $\mathrm{L}_{\theta}$ from $C^{\infty}(M)$ onto $M_{\theta}$ such that 
$\mathrm{L}_{\theta}(f\cdot g)=\mathrm{L}_{\theta}(f) \times_{\theta} \mathrm{L}_{\theta}(g)$. Then it is identified with the fixed point algebra of the set of all 
smooth maps from M to the noncommutative 2-torus $T^2_{\theta}$ under the 
diagonal tensor action $\varphi \otimes \alpha^{-1}$ of $T^2$ where $\alpha$ is the gauge action of $T^2$ on $T^2_{\theta}$. Since $M$ is 4-dimensional, its 
Hodge *-operation on the Grassmann algebra $\Omega(M)$ of all forms of $M$ could be shifted on the all forms $\Omega(M_{\theta})$ of $M_{\theta}$, which is 
denoted by $*_{\theta}$. Let us take a noncommutative principal bundle as 
follows: let 
\[ G \stackrel{\varrho} \rightarrow N \stackrel{\pi} \rightarrow M \]
\noindent
be a principal $G$-bundle over $M$. Suppose there exists a smooth action 
$\tilde{\varphi}$ of a covering group $\tilde{T^2}$ of $T^2$ to 
$\mathrm{Iso}(N,\pi^*(g))$ commuting with $\varrho$, then it follows from [LS] 
that there exists a smooth action $\tilde{\varrho}$ of $G$ on 
$N_{\tilde{\theta}}$ such that
\[G \stackrel{\tilde{\varrho}} \rightarrow N_{\tilde{\theta}}
\stackrel{\pi_{\theta}} \rightarrow M_{\theta}  \]
\noindent
is a noncommutative $G$-bundle over $M_{\theta}$ where $\pi^*(g)$ means the 
pull back Riemannian metric of $g$ on $N$ under $\pi$. Let $\sigma$ be a highestweight of $G$ with respect to its maximal torus $T$ and $V_{\sigma}$ its 
irreducible $G$-module. We identify $\sigma$ with its associated irreducible 
representation of $G$. Let $\Xi_{\sigma}$ be the fixed point algebra of 
$N_{\tilde{\theta}}\otimes V_{\sigma}$ under the diagonal tensor action 
$\tilde{\varrho}\otimes \sigma^{-1}$ of $G$. Then it is a finitely generated 
projective irreducible right $M_{\theta}$-module, and there exists a natural 
number $n$ and a projection $P_{\sigma}$ in $M_n(M_{\theta})$ such that 
$\Xi_{\sigma}=P_{\sigma}(M_{\theta}^n)$. We now take the Grassman connection 
$\nabla_{\sigma}=P_{\sigma}d_{\theta}^n$ of $\Xi_{\sigma}$ where $d_{\theta}$ 
is the canonical outer derivative of $\Omega(M_{\theta})$ derived by 
the ordinary outer derivative $d$ of $\Omega(M)$. We have nothing to know at 
the moment whether $\nabla_{\sigma}$ is self-dual or anti selfdual although this is the case in the special set up (for instance [CDV],[LS]). Let us denote by 
$\mathcal{C}_{+}(\Xi_{\sigma})$ the set of all self-dual compatible connections of $\Xi_{\sigma}$. In what follows, we treat the case that $\mathcal{C}_{+}(\Xi_{\sigma})$ is non empty although we have nothing to check their existence 
for a given principal $G$-bundle over $M$ and an irreducible representation of 
$G$. According to [LS], the set $\mathcal{C}_{\mathrm{YM}}(\Xi_{\sigma})$ of all Yang-Mills connections of $\Xi_{\sigma}$ contains the sum of $\mathcal{C}_{+}(\Xi_{\sigma})$ \\

\Large{\bf{\S3.~Geometry of Noncommutative Instantons}} \large In this section, we analize a geometric structure of $\mathcal{C}_{+}(\Xi_{\sigma})$ defined in 
the previous section under the case where they are non empty. Let $\nabla_{\sigma} \in \mathcal{C}_{+}(\Xi_{\sigma})$ for an irreducible representation 
$\sigma$ of $G$ on a finite dimensional $\mathbb{R}$-vector space $V_{\sigma}$. Let $P_{-}$ be the projection from 
$\Omega^2(M_{\theta},\hat{\Xi}_{\sigma})$ onto 
$\Omega^2_{-}(M_{\theta},\hat{\Xi}_{\sigma})$, where 
$\hat{\Xi}_{\sigma}=\mathrm{End}_{M_{\theta}}(\Xi_{\sigma})$. 
Since $\mathrm{End}_{M_{\theta}}(\Xi_{\sigma})$ is isomorphic to 
$(N_{\tilde{\theta}}\otimes \mathcal{L}(V_{\sigma}))^{\tilde{\varrho}
\otimes Ad(\sigma)(G)}$, 
then its skew-adjoint part $\hat{\Xi}^{sk}_{\sigma}$ is identified with 
$(N^{s}_{\tilde{\theta}}\otimes \mathcal{L}^{sk}(V_{\sigma}))^{\tilde{\varrho}
\otimes Ad(\sigma)(G)}$, 
where $\mathcal{L}(V_{\sigma})$ is the set of all $\mathbb{C}$-linear maps 
on $V_{\sigma}$,~$\mathcal{L}^{sk}(V_{\sigma})$ its skew-adjoint part,~and 
$N^{s}_{\theta}$ is the set of all self-adjoint elements of 
$N_{\tilde{\theta}}$. Since $\mathcal{L}(V_{\sigma})=M_{n_{\sigma}}(\mathbb{C})=\mathcal{U}(n_{\sigma})_{\mathbb{C}}$, then it follows that
$\{T \in \mathcal{L}^{sk}(V_{\sigma})|Tr(T)=0 \}=\mathcal{SU}(n_{\sigma})_{\mathbb{C}}$ where $n_{\sigma}=\mathrm{dim}_{\mathbb{C}}V_{\sigma}$ and $\mathcal{U}(n_{\sigma})_{\mathbb{C}}~,~\mathcal{SU}(n_{\sigma})_{\mathbb{C}}$ are the 
complexifications of the Lie algebras of $\mathrm{U}(n_{\sigma})~,~\mathrm{SU}(n_{\sigma})$ respectively. We then deduce that 
\[\mathcal{U}(n_{\sigma})_{\mathbb{C}}~=~\mathbb{C}\mathrm{I}_{n_{\sigma}}\oplus \mathcal{SU}(n_{\sigma})_{\mathbb{C}} \]
\noindent 
Therefore it implies that 
\[\hat{\Xi}_{\sigma}=M_{\theta}\oplus_{M_{\theta}} (N_{\tilde{\theta}} 
\otimes \mathcal{SU}(n_{\sigma})_{\mathbb{C}})^{\tilde{\varrho}\otimes Ad(\sigma)(G)} \]
\noindent
Let us define
\[\Gamma(ad_{\sigma}(N_{\tilde{\theta}}))=(N_{\tilde{\theta}}\otimes \mathcal{U}(n_{\sigma})_{\mathbb{C}})^{\tilde{\varrho} \otimes Ad(\sigma)(G)}. \]
\noindent
Moreover, we put
\[\Omega^0(ad_{\sigma}(N_{\tilde{\theta}}))=\Gamma(ad_{\sigma}(N_{\tilde{\theta}})),~\Omega^1(ad_{\sigma}(N_{\tilde{\theta}}))=\Omega^1(M_{\theta},\Gamma(ad_{\sigma}(N_{\tilde{\theta}}))) \]
\noindent
,~and 
\[\Omega^2_{-}(ad_{\sigma}(N_{\tilde{\theta}}))=
\mathrm{P}_{-}\Omega^2(M_{\theta},\Gamma(ad_{\sigma}(N_{{\tilde{\theta}}}))).\]
\noindent 
We then introduce an inner product $<~|~>$ on 
$\Omega^j(ad_{\sigma}(N_{\tilde{\theta}}))~(j=0,1,2)$ by 
\[<~\omega~|~\eta~>=\int_{M_{\theta}}\mathrm{tr}(*_{\theta}(\omega^* *_{\theta}\eta)) \]
\noindent
where $*_{\theta}$ is the Hodge operation on $\Omega^*(M_{\theta})$,and $\mathrm{tr}$ is the canonical trace on $M_{n_{\sigma}}(\mathbb{C})$. Using this inner 
product, we also induce the metric topology on $\mathcal{C}_{\pm}(\Xi_{\sigma})$ by the following lemma:\\

\large{\bf{Lemma 3.1}}(cf:[T])~~Let $\nabla_{\sigma} \in 
\mathcal{C}_{+}(\Xi_{\sigma})$. Then it follows that 
\[\mathcal{C}_{+}(\Xi_{\sigma})
=\nabla_{\sigma}+\Omega^1_{+}(ad_{\sigma}(N_{\tilde{\theta}}) \]
\noindent
where ~$\Omega^1_{+}(ad_{\sigma}(N_{\tilde{\theta}}))$ is 
the set of all $\omega \in \Omega^1(ad_{\sigma}(N_{\tilde{\theta}}))$ 
satisfying the equation:
\[  \hat{\nabla}^{-}(\omega)+(\omega^2)_{-}=0 ~,\]
\noindent
and $(\omega^2)_{-}=\mathrm{P}_{-}(\omega^2)$.\\

Proof. Let $\nabla \in \mathcal{C}_{+}(\Xi_{\sigma})$ and put $\omega=\nabla-\nabla_{\sigma}$, then it implies by definition that 
$\omega \in \hat{\Xi}_{\sigma}$. Since $\nabla,~\nabla_{\sigma} \in \mathcal{C}_{+}(\Xi_{\sigma})$, it follows that
\[\mathrm{P}_{-}[\nabla,[\nabla,X]]=\mathrm{P}_{-}[F_{\nabla},X]=[F^{-}_{\nabla},X]=0 \] 
\noindent
for all $X \in \hat{\Xi}_{\sigma}$ respectively, and the similar statement holds for $\nabla_{\sigma}$. We compute that 
\[F_{\nabla}=F_{\nabla_{\sigma}}+[\nabla_{\sigma},\omega]+\omega^2 \], 
\noindent
which implies the conclusion. ~~~ Q.E.D.  \\

As $M$ is compact oriented, it has a spin$^c$ structure. As the similar case as  a spin structure, it follows from [CD$\cdot$V] that there exists a noncommutative spin$^c$ structure $S(c)_{\theta}$ of $M_{\theta}$ as a right 
$M_{\theta}$-module induced by the given spin$^c$ one $S(c)$ of $M$. 
Let $S(c)^{\pm}_{\theta}$ be the half spin structures of $M_{\theta}$ induced by the given half spin$^c$ one $S(c)^{\pm}$ of $M$ respectively. We use the same 
notations as their smooth sections on $M_{\theta}$. 
Then the following statement is easily seen by their definition: \\

\large{\bf{Lemma 3.2}}(cf:[T])~~
\[(1)~~\Omega^1(M_{\theta})=S(c)^{\pm}_{\theta}\otimes S(c)^{\mp}_{\theta}\]
\noindent
\[~~~~~~~~~(2)~~M_{\theta}\oplus \Omega^2_{\pm}(M_{\theta})=
S(c)^{\pm}_{\theta} \otimes S(c)^{\pm}_{\theta} \]
\noindent
respectively.\\

Proof.~(1): ~By definition, $S(c)^{\pm}_{\theta}$ is the algebra of all 
even (odd) polynomials of $\Omega^1(M_{\theta})$ with the Clifford 
multiplication $c_{\theta}$ respectively. 
Then $c_{\theta}(\Omega^1(M_{\theta})S(c)^{\pm}_{\theta}=S(c)^{\mp}_{\theta}$ 
respectively, which implies the conclusion. ~(2):~The result is true if 
$\theta=0$. This isomorphism can also be chosen $\tilde{T}^2$-equivariantly, 
where $\tilde{T}^2$ is a covering group of $T^2$. Then the statement follows. 
~~~Q.E.D.\\

\noindent
By Lemmas 3.1, $\mathcal{C}_{+}(\Xi_{\sigma})$ is identified with 
$\Omega^1_{+}(ad_{\sigma}(N_{\tilde{\theta}}))$ as an affine space respectively. We then show the next three lemma which seems to be quite useful showing our 
main result:\\

\large{\bf{Lemma 3.3}}(cf:[T])~~Let $T_{\nabla}(\mathcal{C}_{+}(\Xi_{\sigma}))$ be the tangent space of $\mathcal{C}_{+}(\Xi_{\sigma})$ at 
$\nabla \in \mathcal{C}_{+}(\Xi_{\sigma})$. Then it follows that 
\[T_{\nabla}(\mathcal{C}_{+}(\Xi_{\sigma}))=\{\omega \in 
\Omega^1(ad_{\sigma}(N_{\tilde{\theta}}))~|~\mathrm{P}_{-}[\nabla,\omega]=0 \}~.\]

Proof.~~Let $\omega_t=\nabla_t-\nabla$ for a smooth curve $\nabla_t \in 
\mathcal{C}_{+}(\Xi_{\sigma})$ with $\nabla_0=\nabla$. 
Then it follows from definition that 
$\omega_t \in \Omega^1(ad_{\sigma}(N_{\tilde{\theta}}))$. We put $\omega$ 
the derivative of $\omega_t$ at $t=0$ in 
$\Omega^1(ad_{\sigma}(N_{\tilde{\theta}}))$. 
Then we see that
\[F_{\nabla_t}=F_{\nabla}+[\nabla,\omega_t]+\omega_t^2 ~,\] 
\noindent
which deduces $F_{\nabla_t}\prime(0)=[\nabla,\omega]$ 
taking their derivatives at $t=0$ since $\omega_0=0$.
As $\nabla_t \in \mathcal{C}_{+}(\Xi_{\sigma})$, it implies that 
$F_{\nabla_t} \in \Omega^2_{+}(ad(N_{\tilde{\theta}}))$, which means that 
$\mathrm{P}_{-}[\nabla,\omega]=0$. ~~~Q.E.D. \\

\noindent
Let $\Gamma(\Xi_{\sigma})$ be the gauge group acting on the set 
$\mathcal{C}(\Xi_{\sigma})$ of all compatible connections of 
$\Xi_{\sigma}$ by 
$\gamma_u(\nabla)=u{\nabla}u^*,(u \in \Gamma(\Xi_{\sigma},\nabla \in \mathcal{C}(\Xi_{\sigma}))$,~and 
$\Gamma(\Xi_{\sigma}) \cdot \nabla$ the orbit of $\nabla$ under $\Gamma(\Xi_{\sigma})$. Then we show the following lemma:\\

\large{\bf{Lemma 3.4}}(cf:[T])~~~~~Let 
$\mathrm{T}_{\nabla}(\Gamma(\Xi_{\sigma})\cdot\nabla)$ be the tangent space 
of $\Gamma(\Xi_{\sigma})\cdot\nabla$ at 
$\nabla \in \mathcal{C}(\Xi_{\sigma})$. 
Then we have that
\[\mathrm{T}_{\nabla}(\Gamma(\Xi_{\sigma})\cdot\nabla)=\{~[\nabla,X]~|~X \in 
\Omega^0(ad(N_{\tilde{\theta}}))~\} \]
\noindent
for any $\nabla \in \mathcal{C}_{+}(\Xi_{\sigma})$.\\

Proof.~~Let $X \in \Omega^0(ad_{\sigma}(N_{\tilde{\theta}}))$ and 
put $\varphi_t=e^{-tX} \in \Gamma(\Xi_{\sigma})$ 
for all $t \in \mathbb{R}$. Then we know that 
\[\gamma_{\varphi_t}(\nabla)(\xi)=e^{-tX}\cdot\nabla(e^{tX}\xi)=t[\nabla,X](\xi)+\nabla(\xi) \]
\noindent
for all $\xi \in \Xi_{\sigma}$. Taking their derivatives at $t=0$, we have that \[ \frac{d}{dt}\gamma_{\varphi_t}(\nabla)\Big{|}_{t=0}(\xi)=[\nabla,X](\xi)\]
\noindent
for all $\xi \in \Xi_{\sigma}$. Therefore the conclusion follows.~~~Q.E.D. \\

\noindent
As we know that 
$\Gamma(\Xi_{\sigma})\cdot\nabla \subseteq \mathcal{C}_{+}(\Xi_{\sigma})$ 
for any $\nabla \in \mathcal{C}_{+}(\Xi_{\sigma})$, we may imagine by Lemma 3,3 and 3.4 the next corollary: \\

\large{\bf{Corollary 3.5}}(cf:[T])~~Let $\mathcal{M}_{\pm}(\Xi_{\sigma})$ be the moduli space of $\mathcal{C}_{+}(\Xi_{\sigma})$ by $\Gamma(\Xi_{\sigma})$. Then the tangent space$T_{[\nabla]}(\mathcal{M}_{+}(\Xi_{\sigma}))$ of $\mathcal{M}_{+}(\Xi_{\sigma})$ at $[\nabla] \in \mathcal{M}_{+}(\Xi_{\sigma})$ is isomorphic to $T_{\nabla}(\mathcal{C}_{+}(\Xi_{\sigma}))/\mathrm{T}_{\nabla}(\Gamma(\Xi_{\sigma})\cdot\nabla)$ as a $\mathbb{C}$-linear space.\\

\noindent
Let us define $\hat{\nabla}^{-}(\omega)=\mathrm{P}_{-}[\nabla,\omega]$ 
for all $\omega \in \Omega^1(ad_{\sigma}(N_{\tilde{\theta}}))$. We then can show the following lemma which is well known in undeformed cases by [AHS] and in a 
special deformed case by [LS]:\\

\large{\bf{Lemma 3.6}}(cf:[T])~~
\[0 \rightarrow \Omega^0(ad_{\sigma}(N_{\tilde{\theta}}))
 \stackrel{\hat{\nabla}} 
\rightarrow \Omega^1(ad_{\sigma}(N_{\tilde{\theta}})) 
\stackrel{\hat{\nabla}^{-}} \rightarrow 
\Omega^2_{-}(ad_{\sigma}(N_{\tilde{\theta}})) \rightarrow 0 \]
\noindent
for $\nabla \in \mathcal{C}_{+}(\Xi_{\sigma})$, where $\hat{\nabla}(X)=[\nabla,X]$ for $X \in \Omega^0(ad_{\sigma}(N_{\tilde{\theta}}))$ and $\hat{\nabla}^{-}(\omega)=\mathrm{P}_{-}[\nabla,\omega]$ for $\omega \in \Omega^1(ad_{\sigma}(N_{\tilde{\theta}}))$. \\

Proof.~~We first show that $\hat{\nabla}^{-}\cdot \hat{\nabla}=0$. Indeed, as 
$\nabla \in \mathcal{C}_{+}(\Xi_{\sigma})$,
\[\hat{\nabla}^{-}\cdot \hat{\nabla}(X)=P_{-}[\nabla,[\nabla,X]]=P_{-}[F_{\nabla},X]=[F^{-}_{\nabla},X]=0 \]
\noindent
for any $X \in \Omega^0(ad_{\sigma}(N_{\tilde{\theta}}))$. We next show that$\mathrm{Ker}~\hat{\nabla}=0$. In fact, suppose $X \in \mathrm{Ker}~\hat{\nabla}$. 
Since $\hat{\nabla}$ commutes with the action $\alpha$ of $T^2$ and $ad(N_{\tilde{\theta}})$ has $\mathcal{U}(n_{\sigma})_{\mathbb{C}}$ as fibres, 
then it follows that 
\[\hat{\nabla}=d_{\tilde{\theta}}^n+\omega_n \] 
\noindent
for some $\omega_n \in M_n(\Omega^1(N_{\tilde{\theta}}))$, 
where $(\omega_n)_{j,k}$ are all central in $\Omega^1(N_{\tilde{\theta}})$, 
which means that $\mathrm{L}_{\theta}(\omega_n)=\omega_n$. Hence it implies that $\hat{\nabla}(X)=\mathrm{P}_{\sigma}\cdot d_{\theta}^n(X)=0$. 
Since $X=\mathrm{L}_{\theta}(X^0)$ for a $X^0 \in \Omega^0(ad(N))$, 
then we have that 
\[ \hat{\nabla}(X)=\mathrm{P}_{\sigma}\cdot d_{\theta}^n(X)=\mathrm{L}_{\theta}
\cdot \mathrm{p}_{\sigma}\cdot d_M^n(X^0)=0~, \]
\noindent
where $\mathrm{p}_{\sigma}$ is the projection in $M_n(C^{\infty}(M))$ 
corresponding to the right $C^{\infty}(M)$-module $\Xi=
(C^{\infty}(N)\otimes V_{\sigma})^{\varrho \otimes \sigma(G)}$ and $d_M$ is 
the outer derivative of $M$. Then $\mathrm{p}_{\sigma}\cdot d_M(X^0)=0$, 
which implies that $X^0=0 \in \Omega^0(ad(N))$ by the irreducibility of 
$\mathrm{p}_{\sigma}\cdot d_M$ on $\Xi$. This means that 
$\mathrm{Ker}~\hat{\nabla}=0$.  ~~~~~Q.E.D.  \\

\noindent
Let us define ~$\mathbb{H}^1_{-}
=\mathrm{Ker}~\hat{\nabla}^{-}/\mathrm{Im}~\hat{\nabla}$~and~
$\mathbb{H}^2_{-}
=\Omega^2_{-}(ad_{\sigma}(N_{\tilde{\theta}})/\mathrm{Im}~\hat{\nabla}^{-}$. 
It follows from Lemma 3.5 that $\mathbb{H}^1_{-}$ is isomorphic to 
$T_{[\nabla]}(\mathcal{M}_{+}(\Xi_{\sigma}))$ as a $\mathbb{C}$-linear space 
for any $\nabla \in \mathcal{C}_{+}(\Xi_{\sigma})$. Moreover, we introduce the 
Laplace type $M_{\theta}$-operators $\Delta^j_{-}$ on 
$\Omega^j(ad_{\sigma}(N_{\tilde{\theta}}))$ \\
\noindent
for $j=1,2$ as follows: 
let $\Delta^{-}_j=(\hat{\nabla}^{-})^*\cdot \hat{\nabla}^{-}+\hat{\nabla}\cdot (\hat{\nabla})^*$ on $\Omega^1(ad_{\sigma}(N_{\tilde{\theta}}))$ and 
$\Delta^{-}_j=\hat{\nabla}^{-}\cdot (\hat{\nabla}^{-})^*$ on 
$\Omega^2(ad_{\sigma}(N_{\tilde{\theta}}))$. 
Then we easily observe the following lemma: \\

\large{\bf{Lemma 3.7}}(cf:[T])~~$\mathbb{H}^j_{-}$ are isomorphic to 
$\mathrm{Ker}~\Delta^{-}_j~(j=1,2)$ as a right $M_{\theta}$-module respectively, and they are finitely generated projective right $M_{\theta}$-modules.\\

Proof.~~$\mathbb{H}^1_{-}$ is $M_{\theta}$-isomorphic to 
$\mathrm{Ker}~\hat{\nabla}^{-} \cap \mathrm{Ker}~(\hat{\nabla})^*$, which is 
equal to $\mathrm{Ker}~\Delta^{-}_j$. The similar way is also valid for 
$\mathbb{H}^2_{-}$. As $\Delta^{-}_j$ are elliptic, the rest is well known.
~~~Q.E.D. \\

\noindent
By Lemma 3.6, the elliptic complex can be described by the following generalized signature $M_{\theta}$-operator:
\[ \hat{\nabla}^{-}+(\hat{\nabla})^*:\Omega^1(ad_{\sigma}(N_{\tilde{\theta}})) \rightarrow \Omega^0(ad_{\sigma}(N_{\tilde{\theta}})) 
\oplus \Omega^2_{-}(ad_{\sigma}(N_{\tilde{\theta}})) \]
\noindent
for all $\nabla \in \mathcal{C}_{+}(\Xi_{\sigma})$. By Lemma 3.7, $\mathbb{H}^j_{-}$ induce the $\mathbb{K}_0(M_{\theta})$-element $[\mathbb{H}^j_{-}]~(j=1,2)$ whose K-theoretic ranks are finite. Let us define the K-theoretic index of $\hat{\nabla}^{-}+(\hat{\nabla})^*$ as follows:
\[\mathrm{Index}_{M_{\theta}}(\hat{\nabla}^{-}+(\hat{\nabla})^*)=[\mathrm{Ker}(\hat{\nabla}^{-}+(\hat{\nabla})^*)]-[\mathrm{Coker}(\hat{\nabla}^{-}+(\hat{\nabla})^*)] \]
\noindent 
Using Lemma 3.6, we then obtain the following lemma which is a noncommutative 
version of Atiyah-Singer index theorem due to Connes [C]~(cf:[MS]): \\

\large{\bf{Lemma 3.8}}(cf:[T])~~
\[ \mathrm{Index}_{M_{\theta}}(\hat{\nabla}^{-}+(\hat{\nabla})^*)
=~[\mathbb{H}^1_{-}]~-~[\mathbb{H}^2_{-}] \]
\noindent
Let $[M_{\theta}]$ be the fundamental cyclic cocycle of $M_{\theta}$ which is 
essentially defined by the JLO-cocycle appeared in [C]. 
Then it follows by Lemma 3.8 that \\

\large{\bf{Corollary 3.9}}~~
\[<[M_{\theta}],\mathrm{Index}_{M_{\theta}}(\hat{\nabla}^{-}+(\hat{\nabla})^*)>
=\mathrm{rank}_{M_{\theta}}[\mathbb{H}^1_{-}]~-~\mathrm{rank}_{M_{\theta}}[\mathbb{H}^2_{-}] \]
\noindent
On the other hand, we consider the Dirac type $M_{\theta}$-operator 
$D^{+}_{\nabla}$ on $\Omega^1(ad_{\sigma}(N_{\tilde{\theta}}))$ 
associated with $\nabla \in \mathcal{C}_{+}(\Xi_{\sigma})$ 
defined by the following process: By Lemma 3.2.(1), \\

$\Omega^1(ad_{\sigma}(N_{\tilde{\theta}})
=\Omega^0(ad_{\sigma}(N_{\tilde{\theta}}))
\otimes (S(c)^{+}_{\theta}\otimes S(c)^{-}_{\theta})$ \\

$\stackrel{\hat{\nabla}} \rightarrow 
\Omega^0(ad_{\sigma}(N_{\tilde{\theta}}))\otimes \Omega^1(M_{\theta})
\otimes (S(c)^{+}_{\theta}\otimes S(c)^{-}_{\theta}) $ \\

$\stackrel{\mathrm{Id} \otimes c^{+}_{\theta}} \rightarrow 
\Omega^0(ad_{\sigma}(N_{\tilde{\theta}}))\otimes 
(S(c)^{-}_{\theta}\otimes S(c)^{-}_{\theta})=\Omega^0(ad_{\sigma}(N_{\tilde{\theta}})) \oplus \Omega^2_{-}(ad_{\sigma}(N_{\tilde{\theta}}))$ \\

\large{\bf{Lemma 3.10}}(cf:[T])~~Given a $\nabla \in \mathcal{C}_{+}(\Xi_{\sigma})$,
\[{\mathrm{Index}}_{M_{\theta}}(D^{+}_{\nabla})={\mathrm{Index}}_{M_{\theta}}(\hat{\nabla}^*+\hat{\nabla}^{-}) \] 
\noindent 
in $K_0(M_{\theta})$. \\ 

Proof.~~Since $\Omega(M_{\theta})\simeq \Omega(M)$ as vector spaces and 
$\hat{\nabla}$ and $*$-operation commute with $T^2$-action $\alpha$, then 
the operator $\hat{\nabla}^*+\hat{\nabla}^{-}$ can be shifted to the elliptic 
operator $(d_M)^*+p_{-}\cdot d_M$ from $\Omega^1(M)$ to 
$\Omega^0(M) \oplus \Omega^2_{-}(M)$. Let $P^0_{-}$ be the projection from 
$\Omega^2(M)$ to $\Omega^2_{-}(M)$ respectively. Using the same argument as 
in [AHS] and taking their principal symbols, it follows from Lemma 3.2 that 
$(d_M)^*+P^0_{\mp}\cdot d_M$ is identified with the Dirac operator $D^{+}$ 
in the following:
\[D^{+}:~S(c)^{+} \otimes S(c)^{-} \rightarrow S(c)^{-} \otimes S(c)^{-} \]~,
\noindent
so far taking their $C^{\infty}(M)$-indices. Then 
$\hat{\nabla}^*+\hat{\nabla}^{-}$ is identified with $D^{+}_{\sigma}$ 
defined before:
\[ D^{+}_{\sigma}: \Omega^0(ad_{\sigma}(N_{\tilde{\theta}}))
\otimes (S(c)^{+}_{\theta} \otimes S(c)^{-}_{\theta}) \rightarrow 
\Omega^0(ad_{\sigma}(N_{\tilde{\theta}}))
\otimes (S(c)^{-}_{\theta} \otimes S(c)^{-}_{\theta}) \]
\noindent
so far taking their $M_{\theta}$-indices. This implies the conclusion.~Q.E.D.\\

\noindent 
In what follows, we want to determine the geometric structure of the moduli 
space $\mathrm{M}_{+}(\Xi_{\sigma})$ of $\mathrm{C}_{+}(\Xi_{\sigma})$ by 
the gauge group $\Gamma(\Xi_{\sigma})$. First of all, we introduce a 
noncommutative Kuranishi map in the following process : 
Given a $\nabla \in \mathcal{C}_{+}(\Xi_{\sigma})$, 
we define a Sobolev norm $||\cdot||_p,~(p\gg 4)$ on 
$\Omega^1_{+}(ad_{\sigma}(N_{\tilde{\theta}})$ by \\

$~~~~~~~~~~~~~~~~~~||\omega||_p = \sum^p_{k=0}\int_{M_{\theta}}~
||(\hat{\nabla})^k(\omega)||^2$ \\

\noindent
for $\omega \in \Omega^1_{+}(ad_{\sigma}(N_{\tilde{\theta}}))$. 
Let $S^p_{+}$  be the completion of 
$\Omega^1_{+}(ad_{\sigma}(N_{\tilde{\theta}}))$ with respect to $||\cdot||_p$. 
By their definition, $\hat{\nabla},~\hat{\nabla}^*$ are bounded from 
$S^p_{+}$ to $S^{p+1}_{+},~S^{p-1}_{+}$ respectively. \\

\noindent 
Since  $\Delta^{-}_j$ is elliptic in $\mathrm{End}_{M_{\theta}}(S^p_{+})$, 
it is a $M_{\theta}$-Fredholm operator on $S^p_{+}$ as well. 
Let $\mathrm{P}_{\mathbb{H}^2_{-}}$ be the projection on $\mathbb{H}^2_{-}$, and define the Green $M_{\theta}$-operator $G^{-}$ of $\Delta^{-}_2$ on $S^p_{+}$ 
as follows :
\[G^{-}=(\mathrm{Id}_{S^p_{+}}-\mathrm{P}_{\mathbb{H}^2_{-}})\cdot 
(\Delta^{-}_2)^{-1}= (\Delta^{-}_2)^{-1} \cdot 
(\mathrm{Id}_{S^p_{+}}-\mathrm{P}_{\mathbb{H}^2_{-}}). \]
\noindent
Then we easily see that \\

$(1):~\mathrm{Im}~G^{-}=\mathrm{Im}~\Delta^{-}_2$, \\

$(2):~G^{-} \cdot \Delta^{-}_2=\Delta^{-}_2 \cdot G^{-}$, \\

$(3):~G^{-} \cdot \Delta^{-}_2=
\mathrm{Id}_{S^p_{+}}-\mathrm{P}_{\mathbb{H}^2_{-}}$ \\

\noindent
We now introduce a densely defined map $\Pi^{-}_{\nabla}$ 
on $S^p_{+}$ by
\[\Pi^{-}_{\nabla}(\omega)= 
\omega+(\hat{\nabla}^{-})^* \cdot G^{\pm}\{(\omega^2)_{-}\} \]
\noindent
for all $\omega \in \Omega^1_{+}(ad_{\sigma}(N_{\tilde{\theta}}))$ 
respectively, which is called a noncommutative Kuranishi map, which is known 
in commutative cases (cf:[AHS]). Then we deduce the following observation: \\

\large{\bf{Lemma 3.11}}(cf:[T])~~\\

$~~~~~~\Omega^1_{+}(ad_{\sigma}(N_{\tilde{\theta}}))=
\mathrm{Ker}~(\hat{\nabla}^{-} \cdot \Pi^{-}_{\nabla}) \cap 
\{P_{\mathbb{H}^2_{-}}((\omega^2)_{-})=0 \}$ ~.\\

Proof.~~By Lemma 3.1, $\omega \in \Omega^1_{+}(ad_{\sigma}(N_{\tilde{\theta}}))$ satisfies 
\[\hat{\nabla}^{-}(\omega)+(\omega^2)_{-}=0 ~.\]
\noindent
Then we compute \\

$\hat{\nabla}^{-}\cdot \Pi^{-}_{\nabla}(\omega)=\hat{\nabla}^{-}(\omega)+\Delta^{-}_2\cdot G^{-}((\omega^2)_{-})$ \\

$~~~~=\hat{\nabla}^{-}(\omega)+(\omega^2)_{-}-P_{\mathbb{H}^2_{-}}((\omega^2)_{-}=-P_{\mathbb{H}^2_{-}}((\omega^2)_{-} $ \\

\noindent
, which implies that 
\[\hat{\nabla}^{-}\cdot \Pi^{-}_{\nabla}(\omega)=-P_{\mathbb{H}^2_{-}}((\omega^2)_{-} ~.\]
\noindent
By the definition of $\mathbb{H}^2_{-}$, it follows that 
\[(\hat{\nabla}^{-})^* \cdot \hat{\nabla}^{-}\cdot \Pi^{-}_{\nabla}(\omega)=0 \]\noindent
, which implies by taking their inner products that 
\[\hat{\nabla}^{-}\cdot \Pi^{-}_{\nabla}(\omega)
=P_{\mathbb{H}^2_{-}}((\omega^2)_{-}=0~.\] 
\noindent
The converse implication is also easily seen. ~Q.E.D. \\

\noindent
We then observe together with Lemma 3.1 the following corollary: \\

\large{\bf{Corollary 3.12}}~~If $\nabla \in \mathcal{C}_{+}(\Xi_{\sigma})$, then it follows that
\[\mathcal{C}_{+}(\Xi_{\sigma})=\nabla+\mathrm{Ker}~(\hat{\nabla}^{-} \cdot \Pi^{-}_{\nabla}) \cap \{P_{\mathbb{H}^2_{-}}((\omega^2)_{-})=0 \} ~.\]
\noindent
It also follows from Lemma 3.11 that \\

\large{\bf{Corollary 3.13}}~~Suppose $\mathbb{H}^2_{\mp}=0$ for a $\nabla \in \mathcal{C}_{+}(\Xi_{\sigma})$, then it follows that
\[\Omega^1_{+}(ad_{\sigma}(N_{\tilde{\theta}}))=
\mathrm{Ker}~(\hat{\nabla}^{-}\cdot \Pi^{-}_{\nabla}) \]
\noindent
Moreover, we obtain by using Lemma 3.11 the following lemma: \\

\large{\bf{Lemma 3.14}}(cf:[T])~~
\[\Pi^{-}_{\nabla}(\Omega^1_{+}(ad_{\sigma}(N_{\tilde{\theta}})))=\mathbb{H}^1_{-} \]
\noindent
for all $\nabla \in \mathcal{C}_{+}(\Xi_{\sigma})$, \\

Proof.~~As $\hat{\nabla}^{-}\cdot \hat{\nabla}=0$ for all $\nabla \in \mathcal{C}_{+}(\Xi_{\sigma})$, then $(\hat{\nabla})^*\cdot (\hat{\nabla}^{-})^*=0$. Let $\eta=\Pi^{-}_{\nabla}(\omega)$ for any $\omega \in \Omega^1_{+}(ad_{\sigma}(N_{\tilde{\theta}}))$. Then we have that
\[(\hat{\nabla})^*(\eta)=(\hat{\nabla})^*(\omega)+\Delta^{-}_2\cdot G^{-}(
\omega^2)_{-}=0 \]~,
\noindent
respectively. By Lemma 3.11, ~~
\[ \hat{\nabla}^{-}(\eta)=P_{\mathbb{H}^2_{-}}((\omega^2)_{-})=0 ~.\]
\noindent
Since $\mathrm{Ker}~(\hat{\nabla})^*=\Omega^1(ad_{\sigma}(N_{\tilde{\theta}}))$, then the conclusion follows. ~~~Q.E.D. \\

\noindent
Using the Sobolev space $S^p_{+}$, we compute the Frechet differentiation 
$\frac{d\Pi^{-}_{\nabla}}{dt}\Big{|}_{t=0}$ of 
$\Pi^{-}_{\nabla}$ on $\Omega^1_{+}(ad_{\sigma}(N_{\tilde{\theta}}))$ in the following lemma: \\

\large{\bf{Lemma 3.15}}~~
\[\frac{d\Pi^{-}_{\nabla}}{dt}\Big{|}_{t=0}=
\mathrm{Id} \]
\noindent
on $\Omega^1_{+}(ad_{\sigma}(N_{\tilde{\theta}}))$ with respect to 
$||\cdot||_p$.\\ 

Proof.~~\\

$\frac{d\Pi^{-}_{\nabla}}{dt}\Big{|}_{t=0}(\omega)
=\displaystyle\lim_{t\to 0}~t^{-1}(\Pi^{-}_{\nabla}(t\omega)-\Pi^{-}_{\nabla}(0))$ \\
$~~~~~~~~~~~~~~~~~~~=\displaystyle\lim_{t\to 0}~(\omega+t(\hat{\nabla}^{-})^* 
\cdot G^{-}((\omega^2)_{-}))=\omega $ \\

\noindent
for any $\omega \in \Omega^1_{+}(ad_{\sigma}(N_{\tilde{\theta}}))$, 
which implies the conclusion. ~~~Q.E.D. \\

\noindent
By the inverse function theorem on Banach spaces, there exists a $\epsilon >0$ 
neighborhood $U(0,\epsilon)$ of $0 \in S^p_{+}$ with respect to $||\cdot||_p$
 on which $(\Pi^{-}_{\nabla})^{-1}$ is diffeomorphic. Using this fact, we next 
show the following lemma: \\

\large{\bf{Lemma 3.16}}~~~~~$U(0,\epsilon) \cap \mathbb{H}^1_{-}$ is 
diffeomorphic to $(\Pi^{-}_{\nabla})^{-1}(U(0,\epsilon))$ \\
\noindent
$~~\cap~\Omega^1_{+}(ad_{\sigma}(N_{\tilde{\theta}}))$ 
under $(\Pi^{-}_{\nabla})^{-1}$.\\

Proof.~~By Lemma 3.14, we know that 
\[\Pi^{-}_{\nabla}(\Omega^1_{+}(ad_{\sigma}(N_{\tilde{\theta}})))=\mathbb{H}^1_{-}~, \]
\noindent
which implies the conclusion. ~~~Q.E.D. \\

\noindent 
Let $P\Gamma(\Xi_{\sigma})=\Gamma(\Xi_{\sigma})/Z(\Gamma(\Xi_{\sigma}))$ be 
the  projective gauge group of $\Gamma(\Xi_{\sigma})$ 
where $Z(\Gamma(\Xi_{\sigma})$ is the center of $\Gamma(\Xi_{\sigma})$. 
Since $\Xi_{\sigma}$ is an irreducible right $M_{\theta}$-module, 
$Z(\Gamma(\Xi_{\sigma})$ is actually the center $ZU(M_{\theta}))$ of 
the unitary group $U(M_{\theta})$ of $M_{\theta}$. We define the gauge action 
$\gamma$ of $P\Gamma(\Xi_{\sigma})$ on $\mathcal{C}_{+}(\Xi_{\sigma})$ by 
$\gamma_{[u]}(\nabla)(\xi)=u \nabla u^*(\xi)$ for any $[u]=uZU(M_{\theta})) \in P\Gamma(\Xi_{\sigma}),~\nabla \in \mathcal{C}_{+}(\Xi_{\sigma})$. We then easily see that: \\

\large{\bf{Lemma 3.17}}~~The action $\gamma$ is effective. \\

Proof.~~The statement is almost clear by its definition. ~~~Q.E.D. \\

\noindent
As the similar way to the commutative cases (cf:[AHS],[FU]), we then show 
the so-called slice theorem in noncommutative setting: \\

\large{\bf{Lemma 3.18}}(cf:[T])~~Let $\epsilon$ be in Lemma 3.16. Then there 
exist \\
$0<\delta \leq \epsilon$,~$U(0,\delta)=\{\omega \in \Omega^1(ad_{\sigma}(N_{\tilde{\theta}}))~|~||\omega||_p <\delta \}$,~and a mapping 
\[ \psi:~\nabla+U(0,\delta) \rightarrow P\Gamma(\Xi_{\sigma}) \] 
\noindent
with the property that \\

$(1):~\hat{\nabla}^*(\gamma_{\psi(\nabla+\omega)}^*(\nabla+\omega)-\nabla)=0$ \\

(2):~~~there exists a $P\Gamma(\Xi_{\sigma})$-equivariant diffeomorphism 
\[~~~\Phi:\nabla+U(0,\delta) \rightarrow U(0,\mathrm{Id}) \subseteq 
\mathrm{Ker}~\hat{\nabla}^* \times P\Gamma(\Xi_{\sigma}).\]
\noindent
defined by 
\[\Phi(\nabla+\omega)=(\gamma_{\psi(\nabla+\omega)}^*(\nabla+\omega)-\nabla,\psi(\nabla+\omega)) \]

Proof.~~Let us consider the next equation defined by 
\[ \mathrm{L}_{\hat{\nabla}}(\omega,[u])=\hat{\nabla}^*(u^*{\nabla}u+u^*{\omega}u)=0 \]
\noindent
for $\nabla+\omega \in \mathcal{C}_{+}(\Xi_{\sigma})$ and $[u] \in P\Gamma(\Xi_{\sigma})$. Then the differential $\delta \mathrm{L}_{\hat{\nabla}}$ of 
$\mathrm{L}_{\hat{\nabla}}$ at $(0,\mathrm{Id})$ is given by 
\[\delta \mathrm{L}_{\hat{\nabla}}:\Omega^1_{+}(ad_{\sigma}(N_{\tilde{\theta}})) \times \Omega^0(ad_{\sigma}(N_{\tilde{\theta}})) \rightarrow \Omega^0(ad_{\sigma}(N_{\tilde{\theta}})) \]
\[\delta \omega \oplus \delta [u] \rightarrow \hat{\nabla}^*([\nabla,\delta u]+\delta \omega) \]~,
\noindent
The partial differential $\delta_2\mathrm{L}_{\hat{\nabla}}$ of $\delta \mathrm{L}_{\hat{\nabla}}$ in the second factor is 
$\hat{\nabla}^* \cdot \hat{\nabla}$, which is a self-adjoint elliptic operator. Since $\mathrm{Ker}~\hat{\nabla}=0$, then it follows by standard elliptic theory $\hat{\nabla}^* \cdot \hat{\nabla}$ is invertible on $S^p_{+}$. By the implicit function theorem in Hilbert space, we obtain a neibourhood $U(\omega,\delta)$ 
with $0<\delta \leq \epsilon$ and a map $\psi:~U(\omega,\delta) \rightarrow 
P\Gamma(\Xi_{\sigma})$ satisfying our conditions. Hence we obtain the 
diffeomorphism $\Phi$ on $\nabla+U(0,\delta)$ cited in this lemma. 
Since we see that
\[\mathrm{L}_{\hat{\nabla}}([v]^*\omega,[v^{-1}u])=\mathrm{L}_{\hat{\nabla}}(\omega,[u])~,\]
\noindent
for all $[u],[v] \in P\Gamma(\Xi_{\sigma}),~\omega \in \Omega^1(ad_{\sigma}(N_{\tilde{\theta}}))$, then $\Phi$ is chosen $P\Gamma(\Xi_{\sigma})$-equivariantly. ~~~Q.E.D.  \\

\noindent
By $(2)$ of the above lemma 3.18, $U(0,\delta)$ can be chosen as $P\Gamma(\Xi_{\sigma})$-invariant. We then prove the following lemma: \\

\large{\bf{Lemma 3.19}}(cf:[T])~~Let $\nabla \in \mathrm{C}_{+}(\Xi_{\sigma})$. Then there exists a $\delta >0$ such that 
\[\nabla_{\sigma}+\{U(0,\delta)\cap \Omega^1_{+}(ad_{\sigma}(N_{\tilde{\theta}}))) \cap \mathrm{Ker}~(\hat{\nabla})^*\} \]
\noindent
is diffeomorphically imbedded in the moduli space 
$\mathcal{M}_{+}(\Xi_{\sigma})$ of $\mathcal{C}_{+}(\Xi_{\sigma})$ by $\Gamma(\Xi_{\sigma})$. \\

Proof.~~By Lemma 3.18, $\nabla_{\sigma}+\{U(0,\delta)\cap \Omega^1(ad_{\sigma}(N_{\tilde{\theta}})\}$ is $P\Gamma(\Xi_{\sigma})$-equivariantly diffeomorphic to  $U(0,\mathrm{Id}) \subseteq \mathrm{Ker}~\hat{\nabla}^* \times P\Gamma(\Xi_{\sigma})$. Then we conclude that 
\[(\nabla_{\sigma}+\{U(0,\delta) \cap \mathrm{Ker}~(\hat{\nabla})^*\}) \cap \mathrm{C}_{+}(\Xi_{\sigma})\]
\noindent
is diffeomorphically embedded in $\mathcal{M}_{+}(\Xi_{\sigma})$. ~~Q.E.D.\\

\noindent
We next show the following lemma (cf:[FU]): \\

\large{\bf{Lemma 3.20}}~~$\mathcal{M}_{+}(\Xi_{\sigma})$ is a Hausdorff space.\\

Proof.~~It suffices to show that $\mathcal{M}(\Xi_{\sigma})=\mathcal{C}(\Xi_{\sigma})/\Gamma(\Xi_{\sigma})$ is Hausdorff. So we need to check that 
\[ \Upsilon=\{(\nabla,\gamma_u(\nabla))~|~\nabla \in \mathcal{C}(\Xi_{\sigma}),
u \in \Gamma(\Xi_{\sigma})~\}  \]
\noindent
is closed in $\mathcal{C}(\Xi_{\sigma}) \times \mathcal{C}(\Xi_{\sigma})$. 
In fact, if $\{(\nabla+\omega_j,\gamma_{u_j}(\nabla+\omega_j))\}_j \subseteq \Upsilon$ is convergent to $(\nabla+\omega,\nabla+\omega\prime)$, then $\omega_j \rightarrow \omega$ and 
\[\gamma_{u_j}(\nabla+\omega_j)=\nabla+u_j\hat{\nabla}(u_j^*)+u_j{\omega}u_j^* 
\rightarrow \nabla+\omega\prime ~.\]
\noindent
 Put $\omega\prime_j=u_j\hat{\nabla}(u_j^*)+u_j{\omega}u_j^* \in \Omega^1(ad_{\sigma}(N_{\tilde{\theta}}))$. Then $\hat{\nabla}(u_j^*) =u_j^*\omega\prime-{\omega}u_j^*$. Since $||u_j\nabla||_k~(0\leq k \leq p-1)$ are bouded by definition 
of the Sobolev norms $||\cdot||$, Hence $\{u_j\}$ is bounded in $S^p(\Omega^0(ad_{\sigma}(N_{\tilde{\theta}})))$. Then it follows by Rellich's theorem that 
there exists a subsequence $\{u_k\}$ of $\{u_j\}$ such that 
$\{u_k\}$ is convergent in $S^{p-1}(\Omega^1(ad_{\sigma}(N_{\tilde{\theta}}))$. Hence $u_k^*\omega\prime-{\omega}u_k^*$ converges in $S^{p-1}(\Omega^1(ad_{\sigma}(N_{\tilde{\theta}})))$, so is $\{\hat{\nabla}(u_k)\}_k$ in $S^{p-1}(\Omega^1(ad_{\sigma}(N_{\tilde{\theta}})))$ as well. Therefore $\{u_k\}_k$ converges to 
$u \in S^p(\Omega^0(ad_{\sigma}(N_{\tilde{\theta}})))$. Since $\Gamma(\Xi_{\sigma})$ is closed in $S^p(\Omega^0(ad_{\sigma}(N_{\tilde{\theta}})))$, we conclude that $\gamma_u(\nabla+\omega)=\nabla+\omega\prime$. This completes the proof.~~Q.E.D. \\

We now state our main theorem which is a generalization of the one appeared in 
[T] as follows:\\

\large{\bf{Theorem 3.21}}~~Let $M$ be a compact oriented toric 4-dimensional 
Riemannian manifold, and $G \stackrel{\varrho} \rightarrow P \rightarrow M$ 
a principal $G$-bundle over $M$.where $G$ is a compact connected Lie group. 
Suppose $G \stackrel{\tilde{\varrho}} \rightarrow P_{\tilde{\theta}} \rightarrow M_{\theta}$ exists, then given the finite generated projective irreducible right $M_{\theta}$-module $\Xi_{\sigma}=(P_{\theta}\otimes V_{\sigma})^{\tilde{\varrho}\otimes Ad(\sigma)(G)}$ for a highest weight $\sigma$ of $G$,~
the moduli space $\mathcal{M}_{+}(\Xi_{\sigma})$ of $\mathcal{C}_{+}(\Xi_{\sigma})$ under $\Gamma(\Xi_{\sigma})$ is a locally smooth manifold with its dimension:
\[<[M_{\theta}],\mathrm{ch}_{\theta}(\Xi_{\sigma} \otimes S(c)^{-})> 
+ \mathrm{rank}_{M_{\theta}}[\mathbb{H}^2_{-}] \]
\noindent
where $\mathrm{ch}_{\theta}$ is the Connes-Chern character from \\ $K_0(M_{\theta})$ to the perodic cyclic homology of $M_{\theta}$, and $S(c)^{-}$ is the 
smooth sections of the (-)-half spin$^c$-structure of $M$.\\

Proof. ~~By~Corollary 3.9,~Lemma 3.17,~3.19 and 3.20,~$\mathcal{M}_{+}(\Xi_{\sigma})$ is a smooth manifold.  By Lemma 3.16, its dimension is equal to 
$\mathrm{rank}_{M_{\theta}}\mathbb{H}^1_{-}$, which is by Corollary 3.9 and 
Lemma 3.10 that 
\[<[M_{\theta}],\mathrm{ch}_{\theta}(\Xi_{\sigma} \otimes S(c)^{-})>+\mathrm{rank}_{M_{\theta}}[\mathbb{H}^2_{-}]~.\]
\noindent 
This completes the proof. ~~~Q.E.D. \\

noindent
From this theorem, we deduce several useful corollaries which are well known in  the case of ordinary manifolds as well as a noncommutative 4-sphere case:\\

\large{\bf{Corollary 3.22}}~~In the above theorem, suppose $[\mathbb{H}^2_{-}]=0 \in K_0(M_{\theta})$ for all $\nabla \in \mathcal{C}_{+}(\Xi_{\sigma})$, then 
the instanton moduli space $\mathcal{M}_{+}(\Xi_{\sigma})$ of $Xi_{\sigma}$ is 
a smooth manifold with its dimension:
\[<[M_{\theta}],\mathrm{ch}_{\theta}(\Xi_{\sigma} \otimes S(c)^{-})> ~.\]

Proof. ~~By Corollary 3.12 and Theorem 3.21, the conclusion follows. ~~~Q.E.D.\\

\large{\bf{Corollary 3.23}}~~Let $M$ be a compact oriented toric Riemannian 
4-manifold with positive scalar curvature,and G a compact connected Lie group. 
Suppose there exists a noncommutative principal $G$-bundle 
\[ G \stackrel{\tilde{\varrho}}\rightarrow N_{\tilde{\theta}} \rightarrow M_{\theta} ~,\]
\noindent
then given a highest weight $\sigma$ of $G$ and consider the right $M_{\theta}$-module $\Xi_{\sigma}=(N_{\tilde{\theta}} \otimes V_{\sigma})^{\tilde{\varrho}
\otimes Ad(\sigma)(G)}$, then $\mathcal{M}_{+}(\Xi_{\theta})$ is a smooth 
manifold with dimension:
\[ \mathrm{dim}\mathcal{M}_{+}(\Xi_{\sigma})=<[M_{\theta}],ch_{\theta}(\Xi_{\theta} \otimes S(c)^{-}> \]

Proof.~~By the same method as in [AHS], as the manifold $M$ has a positive 
scalar curvature, we have that $(\mathbb{H}^0)^2_{-}=0$ for the 2-cohomology 
$(\mathbb{H}^0)^2_{-}$ with respect to the Laplacian associated with 
$(d_M)^*+P^0_{-}\cdot d_M$. Since the quantization map $L_{\theta}$ from 
$C^{\infty}(M)$ to $M_{\theta}$ can be lifted to a K-theoretic isomorphism 
(cf:[R]) and $L_{\theta} \cdot d_M = d_{\theta} \cdot L_{\theta}$, then it 
implies that $\mathbb{H}^2_{-}=0$ for all $\nabla \in \mathcal{C}_{+}(\Xi_{\sigma})$. Hence the conclusion follows from Corollary 3.22. ~~~Q.E.D. \\

\noindent
In what follows, we compute the typical example using the main theorem obtained  in \S3 :\\

\Large{\bf{\S4.~Example}} \large In what follows, we compute the typical example using the main theorem obtained in \S3. Let us consider the Hoph bundle:
\[ \mathrm{U}(1) \longrightarrow S^5 \longrightarrow \mathbb{C}\mathrm{P}^2~, \]
\noindent
and its associated noncommutative one:
\[ \mathrm{U}(1) \stackrel{\varrho}\longrightarrow S^5_{\theta} \longrightarrow \mathbb{C}\mathrm{P}^2_{\theta} \]
\noindent
Let $1$ be the highest weight of $U(1)$ and $\Xi_1=(S^5_{\theta} \otimes \mathbb{C})^{\varrho \otimes \sigma_1^{-1}}(U(1))$ be its associated right $\mathbb{C}P^2_{\theta}$-module where $\sigma_1$ is the character of $U(1)$ corresponding to $1$. 
Let $z_j~(j=1,2,3)$ be the 3-generators of $S^5_{\theta}$ and put $P_1=(z_1^*,z_2^*,z_3^*)^t(z_1,z_2,z_3)$ be the projection of $M_3(\mathbb{C}\mathrm{P}^2_{\theta})$ since $(z_1,z_2,z_3)(z_1^*,z_2^*,z_3^*)^t=I$. Then $\Xi_1=P_1(\mathbb{C}\mathrm{P}^2_{\theta})^3$ up to isomorphism. Let $\nabla=P_1\cdot d_{\theta}^3$ 
the Grassmann connection of $\Xi_1$. Then $\nabla \in \mathcal{C}_{+}(\Xi_1)$. Now we compute the Connes-Chern character $ch_{\theta}(\Xi_1)$ of $\Xi_1$. Then 
it follows from [C] (cf:[M]) that\\

$~~~~ch_0(\Xi_1)=Tr(P_1)=\sum_j~z_j^*z_j=I \in \mathbb{C}\mathrm{P}^2_{\theta}~,$ \\

$~~~~ch_1(\Xi)=Tr((P-\frac{1}{2}I)\otimes P^{\otimes 2})~,~ch_2(\Xi_1)=Tr((P-\frac{1}{2}I)\otimes P^{\otimes 4})$ \\

\noindent
where $Tr$ is the canonical trace of $M_3(\mathbb{C})$. We also compute the 
chern character $ch(S(c)^{-})$ of the spin$^c$ bundle $S(c)^{-}$ over $\mathbb{C}\mathrm{P}^2$. Let us take the real Fubini-Study metric as follows:
\[ ds^2=\frac{dr^2+r^2e_z^2}{(1+r^2)^2}+\frac{r^2(e_x^2+e_y^2)}{(1+r^2)} \]
\noindent
where $\rho_1=r\mathrm{cos}(\psi)\mathrm{exp}(i\alpha),\rho_2=r\mathrm{sin}(\psi)\mathrm{exp}(i\beta)$ for any $(\rho_1,\rho_2,1) \in \mathbb{C}\mathrm{P}^2~,$ \\

\noindent
and \\

$e_x=-\mathrm{sin}(\phi)\mathrm{cos}(\phi)\mathrm{cos}(\alpha+\beta)d\alpha+\mathrm{sin}(\phi)\mathrm{cos}(\phi)\mathrm{cos}(\alpha+\beta)d\beta$ \\

$~~~~~~+\mathrm{sin}(\alpha+\beta)d\psi~,$ \\

$e_y=-\mathrm{sin}(\phi)\mathrm{cos}(\phi)\mathrm{cos}(\alpha+\beta)d\alpha
+\mathrm{sin}(\phi)\mathrm{cos}(\phi)\mathrm{cos}(\alpha+\beta)d\beta$ \\

$~~~~~~-\mathrm{cos}(\alpha+\beta)d\psi~,$ \\

$e_z=\mathrm{cos}^2(\phi)d\alpha+\mathrm{sin}^2(\psi)d\beta ~.$ \\

\noindent
We then define the spin$^c$ connection $\nabla$ on the smooth sections of $S(c)^{-}$ as follows: 
\[\nabla = (-\frac{1}{4}\sum_{j,k=0,1,2,3}~\omega_{jk}\Gamma^{jk}+iI) \]
\noindent
where 
\[\Gamma^0=\pmatrix{ I_2 & 0 \cr 0 & -I_2 \cr } ,~~\Gamma^j=\pmatrix{ 0 & 
i\sigma^j \cr i\sigma^j & 0 },~~\Gamma^{jk}= \frac{1}{2}[\Gamma^j,\Gamma^k]. \]
\noindent
for the Pauli matrices $\sigma^j~(j=1,2,3,4)$,~and \\

$\omega_{00}=0~,~\omega_{01}=\frac{-1}{(1+r^2)^{1/2}}e_x~,~\omega_{02}=\frac{-1}{(1+r^2)^{1/2}}e_y~,~\omega_{03}=\frac{r^2-1}{r^2+1}e_z~,$ \\

$\omega_{10}=\frac{-1}{(1+r^2)^{1/2}}e_x~,~\omega_{11}=0~,~\omega_{12}=\frac{1+2r^2}{1+r^2}e_z~,~\omega_{13}=\frac{-1}{(1+r^2)^{1/2}}e_y~,$ \\

$\omega_{20}=\frac{-1}{(1+r^2)^{1/2}}e_y~,~\omega_{21}=\frac{1+2r^2}{1+r^2}e_z~,~\omega_{22}=0~,~\omega_{23}=\frac{1}{(1+r^2)^{1/2}}e_x~,$ \\

$\omega_{30}=\frac{r^2-1}{r^2+1}e_z~,~\omega_{31}=\frac{-1}{(1+r^2)^{1/2}}e_y~,~\omega_{32}=\frac{1}{(1+r^2)^{1/2}}e_x~,~\omega_{33}=0~.$ \\

\noindent
Then we compute the curvature $F_{\nabla}$ of $\nabla$ as follows:\\

$F_{\nabla}=\nabla^2=\frac{1}{8}(\sum_{k=1,2,3}~\omega_{0k}\omega_{0k+1}\Gamma^{0k}\Gamma^{0k+1})+\omega_{0k}\omega_{(k+2)k}\Gamma^{0k}\Gamma^{(k+2)k}$ \\

$~~~~~~+\sum_{j=0,1,2,3}~\omega_{j(j+1)}\omega_{(j+1)(j+2)}\Gamma^{j(j+1)}\Gamma^{(j+1)(j+2)}+\omega_{0
1}\omega_{31}\Gamma^{03}\Gamma^{31})$ \\

\noindent
(mod.3). Then it follows by their definitions that 
\[ ch_0(S(c)^{-})=4~,~Tr(F_{\nabla})=0~,~ch_1(S(c)^{-})=0 \]
\noindent
Moreover, we deduce that 
\[ Tr(F_{\nabla}^2)=0~,~ch_2(S(c)^{-})=0 \]
\noindent
Hence we have that $ch(S(c)^{-})=4$. Then we compute by definition (cf:[C],[G$\cdot$B],[M]) that
\[<[CP^2_{\theta}]~,~ch_{\theta}(\Xi_1 \otimes S(c)^{-})>=4\{1-\frac{1}{2}\varphi^{1}(1,1,1)-6\varphi^{2}(1,1,1,1,1)\}=4 \]
\noindent
where $\varphi^{k}$ is the JLO-cyclic $k$-cocycles ~$(k=1,2)$ of $\mathbb{C}\mathrm{P}^2_{\theta}$. We next compute $\mathbb{H}^2_{-}$. First of all, we know 
that
\[\Gamma(ad_{\sigma_1}(S^5_{\theta}))=(S^5_{\theta} \otimes \mathbb{C})^{\varrho \otimes Ad(\sigma_1)(U(1))}=\mathbb{C}\mathrm{P}^2_{\theta} \]
\noindent
We also check the following:
\[\Omega^2_{-}(ad_{\sigma_1}(S^5_{\theta}))=\Omega^2_{-}(\mathbb{C}\mathrm{P}^2_{\theta})=(\Omega^2_{-}(\mathbb{C}\mathrm{P}^2)\otimes T^2_{\theta})^{\phi_* \otimes \alpha^{-1}(T^2)} \]
\noindent
,and 
\[\mathrm{Im}~\hat{\nabla}_{-}=(\mathrm{Im}~d^{-} \otimes T^2_{\theta})^{\phi_*
\otimes \alpha^{-1}(T^2)}~. \]
\noindent
We then check that 
\[\mathbb{H}^2_{-}=\Omega^2_{-}(ad_{\sigma_1}(S^5_{\theta}))/\mathrm{Im}~\hat{\nabla}^{-}=((\Omega^2_{-}(\mathbb{C}\mathrm{P}^2)/\mathrm{Im}~d^{-})\otimes T^2_{\theta})^{\phi_* \otimes \alpha^{-1}(T^2)} ~,\]
\noindent
which is equal to $(\mathbb{H}^2_{-}(\mathbb{C}\mathrm{P}^2)\otimes T^2_{\theta})^{\phi_* \otimes \alpha^{-1}(T^2)}$. Since $\mathbb{C}\mathrm{P}^2$ has a 
positive scalar curvature, it follows by the Bochner-Weitzenbock formula that 
$\mathbb{H}^2_{-}(\mathbb{C}\mathrm{P}^2)=0$. Therefore, we conclude that 
$\mathbb{H}^2_{-}=0$. We then obtain the following result:\\

\large{\bf{Theorem 4.1}}~~Let us take the following Hopf bundle:
\[ \mathrm{U}(1) \longrightarrow S^5 \longrightarrow \mathbb{C}\mathrm{P}^2~,\]
\noindent
and its associated noncommutative one:
\[ \mathrm{U}(1) \stackrel{\varrho}\longrightarrow S^5_{\theta} \longrightarrow  \mathbb{C}\mathrm{P}^2_{\theta} \]
\noindent
Let $1$ be the highest weight of $U(1)$ and $\Xi_1=(S^5_{\theta} \otimes \mathbb{C})^{\varrho \otimes \sigma_1^{-1}}(U(1))$ be its associated right $\mathbb{C}\mathrm{P}^2_{\theta}$-module where $\sigma_1$ is the character of $U(1)$ 
corresponding to $1$. Then its instanton moduli space $\mathcal{M}_{+}(\Xi_1)$ 
is a smooth 4-manifold.\\

\large{\bf{Remark}}~~For the highest weight $n$ of $U(1)$, let $\Xi_n$ be its 
corresponding right $\mathbb{C}\mathrm{P}^2_{\theta}$-module. Then the instanton moduli space $\mathcal{M}_{+}(\Xi_n)$ is a smooth 4n-manifold using the similar method as for $\Xi_1$. \\

\vspace{5mm}

\begin{center}
\Large{References}
\end{center}
\vspace{2mm}
\noindent
[AHS]~M.F.Atiyah, N.J.Hitchin and I.M.Singer,~Self duality in 4-\\
~~~~~~~~dimensional Riemannian geometry, Proc.R.S.London A362 \\
~~~~~~~~(1978), 425-461. \\

\noindent
[C]~A.Connes,~Noncommutative Geometry, Academic Press (1994). \\

\noindent
[CD$\cdot$V]~Connes and M.Dubois$\cdot$Violette,~Noncommutative Finite- \\
~~~~~~~~~Dimensional Manifolds. Spherical Manifolds and Related \\
~~~~~~~~~Examples, Comm.Math.Phys., 230 (2002), 539-579. \\

\noindent
[CL]~A.Connes and G.Landi,~Noncommutative manifolds: The \\
~~~~~~instanton algebra and isospectral deformations, Comm.Math.\\
~~~~~~Phys., 221 (2001), 141-159. \\

\noindent
[FU]~D.S.Freed and K.K.Uhlenbeck,~Instanton and 4-Manifolds, \\
~~~~~~MSRI.Publ.1, Springer-Verlag (1984). \\

\noindent
[G$\cdot$BVF]~J.M.Gracia$\cdot$Bondia, J.C.Varilly and H.Figueroa,~Elements \\
~~~~of Noncommutative Geometry, Birkhauser Advanced Texts (2001). \\

\noindent
[LS]~G.Landi and W.van Suijlekom,~Noncommutative instantons \\
~~~~~~from twisted conformal symmetries, arXiv.math QA/0601554 \\
~~~~~~(2006). \\

\noindent
[MS]~C.C.Moore and C.Schochet,~Grobal Analysis on Foliated \\
~~~~~~Spaces, MSRI.Publ,.9 Springer-Verlag (1988). \\

\noindent
[M]~J.W.Morgan,~The Seiberg-Witten equations and applications to \\
~~~~~topology of smooth four-manifolds, Princeton University Press \\
~~~~~(1996).\\

\noindent
[T]~H.Takai,~Moduli Spaces of Instantons on Noncommutative 4-\\
~~~~Manifolds, arXiv math.DG/0610536 (2006).

\end{document}